\documentclass[oneside,11pt]{article}
\usepackage[english]{babel}
\usepackage{verbatim}
\usepackage{graphicx}
\usepackage{syntonly}
\usepackage{amsmath}
\usepackage{amsfonts}
\usepackage{amssymb}
\usepackage{amsthm}
\usepackage{mathrsfs}
\usepackage[all]{xy}

\makeatletter

\theoremstyle{plain}
\newtheorem{theorem}{Theorem}[section]
\newtheorem{corollary}[theorem]{Corollary}

\newtheorem{proposition}[theorem]{Proposition}
\theoremstyle{definition}
\newtheorem{definition}[theorem]{Definition}
\newtheorem{example}[theorem]{Example}
\newtheorem{remark}[theorem]{Remark}

\newcommand{\eqnsection}{
\renewcommand{\theequation}{{\thesection.\arabic{equation}}} 
\renewcommand{\con}{\mathrm{conv}}

\makeatletter \csname @addtoreset\endcsname{equation}{section}
\makeatother} 

\begin{document}
\thispagestyle{empty}
\date{}

\title{\textbf{\large ALGEBRAIC AND GEOMETRIC ASPECTS\\[1mm] OF BIPARTITE PLANAR GRAPHS\\[3mm]}}

\author{\textsc{Maurizio Imbesi}\\[2mm]
\small{\em Department of Mathematical and Computer Sciences, Physical and Earth Sciences}\\
\small{\em University of Messina, Viale F. Stagno d'Alcontres 31, I-98166 Messina, Italy}\\
\small{{\tt e-mail address:} maurizio.imbesi@unime.it}\\[6mm]
\textsc{Monica La\,Barbiera}\\[2mm]
\small{\em Department of Electrical, Electronic and Computer Engineering}\\
\small{\em University of Catania, Viale A. Doria 6, I-95125 Catania, Italy}\\
\small{{\tt e-mail address:} monica.labarbiera@unict.it}
}

\maketitle


\vspace{-6mm}
\begin{center}
\rule{4cm}{.3pt}\\[1mm]
\emph{\small Dedicated to the memory of Tania Restuccia}
\end{center}
\vspace{0mm}
\begin{quote}
\textsc{Abstract.} {\small Let $B_{2t}$ be a bipartite planar graph
with an even number of regions. We are able to find bounds for the graded
Betti numbers and the projective dimension of the quotient ring
associated to the graph. We also will investigate the minimal vertex
covers and the maximum matchings related to such a graph. }
\end{quote}
\vspace{1mm}
{\small 2020 \emph{Mathematics Subject Classification}. Primary: 05C99; Secondary: 13D02,13F20.}\\
{\small \emph{Key words and phrases}. Planar graph; standard algebraic invariants; vertex cover.}
\smallskip
\section*{Introduction}
\par
Let $G$ be a graph on a vertex set $\{v_1,\ldots, v_n \}$ and $R =
K[X_1, \ldots , X_n]$ be the polynomial ring over a field $K$, with
variables $X_i$ associated to vertices $v_i$ of $G$. The monomial
ideal $I_G$ of $R$ generated by $\{X_iX_j \ | \{v_i,v_j \} \ \
\textrm{is an edge of} \ \ G \}$ is said the edge ideal of $G$.\\
In this paper we are interested to extract specific information
about some invariants linked to the minimal graded resolution of
$R/I_G$ when $G$ is the bipartite planar graph $B_{2t}$, where
$t\geq 1$ is an integer and $r=2t$ the number of its regions. In
\cite{DG}, the $K$-algebra $K[B_{2t}]$ of the graph
$B_{2t}$ was studied using its geometry and the Hilbert function of
$K[B_{2t}]$ of it was computed.\\ The paper is structured as follows.\\
In Section 1 some preliminary notions about
the planar graphs $B_{2t}$ are given.\\
In Section 2 we study all the graded Betti numbers that appear in the minimal graded
resolution of $R/\mathcal{I}$, where $\mathcal{I}$ is the edge ideal
of $B_{2t}$, using some geometric properties of $B_{2t}$. The graded
Betti numbers determine the rank of the free modules in the minimal
graded resolution of $R/\mathcal{I}$ and in general it is not
possible to give a generic formula to compute them. We are able to
give bounds for them in terms of the number of the regions of
$B_{2t}$. As a particular case we study the second Betti number of
degree $3$ of $R/\mathcal{I}$ linked to the number of the regions
of these planar graphs and give an explicit formula to compute it.\\
The last two sections of the work are devoted to study important
sets associated to a bipartite planar graph $G$: the minimal vertex
covers and the maximum matchings of $G$. The problem of the vertex
covering was intensively studied in \cite{{IL1}, {IL2}, {IL3},
{IL4}, {LB}}. It consists in finding a vertex cover of minimum
cardinality, that is a minimal subset $\mathcal{A}$ of the vertex
set of $G$ such that each edge of $G$ is incident with one vertex in
$\mathcal{A}$. More precisely, we describe the minimal vertex covers
of the bipartite planar graphs $B_{2t}$ and connect to the vertex
covers of $B_{2t}$ some algebraic aspects such as dimension and
height. There is a correspondence among the minimal vertex covers
and the minimal primes of the edge ideal. If all minimal vertex
covers have the same size, then the graph is unmixed. The unmixed
bipartite graphs were characterized in \cite{Vil}, and some
generalizations of them were given in \cite{IL}. We will verify that
the planar graphs $B_{2t}$ are not unmixed. Furthermore, as
algebraic topic, we will compute the dimension of $R/\mathcal{I}$
and establish bounds for the projective dimension of $R/\mathcal{I}$
by connecting the geometry of $B_{2t}$ with graph theoretical
properties. The problem to find maximum matchings for the bipartite
graphs $B_{2t}$ is that of associating the geometry of $B_{2t}$ to
the minimal vertex covers. We prove that the graphs $B_{2t}$, for
$t$ odd, have perfect matchings of K\"{o}nig type (\cite{V}), say a
collection $e_1, \ldots, e_g$ of pairwise disjoint edges such that
the union of the vertices in which $e_1, \ldots, e_g$ are incident
is the vertex set of the graph and $g$ is the height of its edge
ideal. Finally we give a complete description of these matchings.

\section{Preliminary notions}
Let $G$ be a graph with vertex set $V(G) =\{v_1,\ldots,
v_n\}$ and edge set $E(G)$ which consists of pairs $\{v_i,v_j\}$
of adjacent vertices, for some $v_i, v_j \in V(G)$.\\
A graph $G$ on vertices $v_1,\ldots, v_n$ is \textit{complete} if
there exists an edge for all pair $\{ v_i,v_j \}$ of
vertices of $G$. It is denoted by $K_n$.\\
A graph $G$ is \textit{bipartite} if its vertex set $V(G)$ can be
partitioned into disjoint subsets $V_{1}=\{x_{1},\ldots,x_{n}\}$ and
$V_{2}=\{y_{1},\ldots,y_{m}\}$, and any edge joins a vertex
of $V_{1}$ to a vertex of $V_{2}$.\\
A graph $G$ is \textit{complete bipartite} if all the vertices
of $V_{1}$ are joined to all the vertices of $V_{2}$. It is
denoted by $K_{n,m}$.
\begin{definition}
A graph $G$ is said \textit{planar} if it has an embedding in the plane
such that each pair of edges is intersected only in the common vertices.
\end{definition}
\noindent A planar graph is subdivided by its edges into plane regions.
\begin{theorem} [\cite{H}, Theorem 11.13] \label{t1}
A graph is planar if and only if it has no subgraphs containing
$K_5$ and $K_{3,3}$.
\end{theorem}
\noindent The complete graphs $K_5$ and $K_{3,3}$ are the minimal, not planar,
graphs. In fact, it is not possible to represent these graphs in the
plane so that the edges are not intersected only in the vertices.\\
Now we consider the class of bipartite planar graphs $B_{2t}$
introduced in \cite{DG}.\\
Let $B_{2t}$ be the planar graph with $r=2t$ regions,
$t\geq 1$ an integer, on vertex set $V(B_{2t})=\{v_1, \ldots,
v_{3t+3} \}$ and edge set $E(B_{2t})= \{ \{v_i,v_{i+1} \} | 1 \leq i
\leq 3t+2, \ \ i\neq t+1,
2t+2 \}  \cup  \{ \{v_i,v_{i+t+1}\} | 1 \leq i \leq 2t+2 \}$.\\
$B_{2t}$ is a planar graph by Theorem \ref{t1}, for all $t\geq 1$.\\
\noindent $B_{2t}$ is a bipartite planar graph. In fact, the vertex
set of $B_{2t}$ can be partitioned into disjoint subsets $V_{1}$ and
$V_{2}$, with  $N = |V_1|+|V_2|=3t+3$ and $|V_i|$ denotes the number of vertices of $V_i$ for $i=1,2$.\\
We have two cases:\\
$1)$ If $t$ is even and $N=3t+3$ is odd one has  $V_{1}=\{v_{i} \ | \ i \
{\rm odd} ,\, 1 \leq i \leq 3t+3 \}$  with $|V_1|=\frac{3t+4}{2}$ and
$V_{2}=\{v_{i} \ | \ i \ {\rm even},\, 1 \leq i \leq 3t+3 \}$  with
$|V_2|=\frac{3t+2}{2}.$\\
$2)$ If $t$ is odd and $N=3t+3$ is even one has
$V_{1}=\{v_{1},v_{3}, \dots, v_{t} \} \cup
\{v_{2+(t+1)},v_{4+(t+1)}, \dots, v_{t+1+(t+1)} \} \cup
\{v_{1+(2t+2)},v_{3+(2t+2)}, \dots, v_{t +(2t+2)} \}$ and
$V_{2}=\{v_{2},v_{4}, \dots, v_{t+1} \} \cup
\{v_{1+(t+1)},v_{3+(t+1)}, \dots, v_{t+(t+1)} \} \cup
\{v_{2+(2t+2)},v_{4+(2t+2)},$  $ \dots,$ $v_{t+1 +(2t+2)} \}$.\\
Note that $|\{v_{1},v_{3},v_{5},\dots, v_{t} \}|=
|\{v_{2},v_{4},v_{6},\dots, v_{t+1} \}|= \frac{t+1}{2}$, hence one
has $|V_1|=|V_2|= \frac{3t+3}{2}$.\\
Then the graph $B_{2t}$ has vertex set $V(B_{2t})= V_1 \cup V_2$,
with $V_1 \cap V_2= \emptyset$, such that its edges join the vertices
of $V_1$ to vertices of $V_2$ only, as follows by definition of
$E(B_{2t})$.
\begin{example}
$G=B_6$\,, with $V(B_{6})=\{v_1, \ldots, v_{12} \}$ and
$E(B_{6})=\{\{v_1, v_2\},$ $\{v_2, v_3\}, \{v_3, v_4\}, \{v_5, v_6\},
\{v_6, v_7\}, \{v_7, v_8\} ,\{v_9, v_{10}\}, \{v_{10}, v_{11}\}, \{v_{11}, v_{12}\},\\
\{v_{1}, v_{5}\},\{v_{2}, v_{6}\}, \{v_{3}, v_{7}\},
  \{v_{4}, v_{8}\},  \{v_{5}, v_{9}\},  \{v_{6}, v_{10}\},  \{v_{7}, v_{11}\},
  \{v_{8}, v_{12}\} \}.$

\begin{figure}[htbp]
\begin{center}
   \includegraphics[scale=.85]{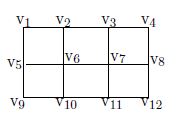}
\end{center}
\end{figure}

\vspace{-8mm}\noindent $V(B_{6})$ can be partitioned into disjoint subsets:\\
$V(B_{6})=\{v_1,v_3,v_6,v_8,v_9,v_{11} \} \cup
\{v_2,v_4,v_5,v_7,v_{10},v_{12} \}= V_1 \cup V_2$.\\
If we rename $\{x_1, \ldots,x_6 \}$ the vertices of  $V_1$  and
$\{y_1, \ldots, y_6 \}$ the vertices of  $V_2$,
then the edge set can be written as:\\[2mm]
$E(B_{6})=\{\{x_1, y_1\}, \{x_2, y_1\}, \{x_2, y_2\}, \{x_3, y_3\},
\{x_4, y_4\}, \{x_5, y_5 \}, \{x_6, y_5 \}, \{x_6, y_6 \},$\\
$\{x_1, y_3 \}, \{x_3, y_1 \}, \{x_2, y_4 \}, \{x_4, y_2 \}, \{x_5,
y_3 \}, \{x_3, y_5 \}, \{x_6, y_4 \}, \{x_4, y_6 \}, \{x_3, y_4
\}\}.$

\smallskip
\begin{figure}[htbp]
\begin{center}
   \includegraphics[scale=.85]{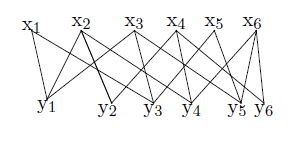}
\end{center}
\end{figure}

\vspace{-8mm}\noindent The two pictures represent the same planar graph $B_6$.
\end{example}
\noindent Let $R = K[X_1, \ldots , X_n]$ be the polynomial ring over a field
$K$, with one variable $X_i$ for each vertex $v_i$ of $G$.
\begin{definition}
We call \textit{edge ideal $I_G$} associated to a graph $G$ the ideal
of $R$ generated by monomials of degree two, $X_iX_j$, on the
variables $X_1, \ldots, X_n$, such that $\{v_i,v_j\} \in E(G)$, for
$1\leq i, j\leq n$.
\end{definition}
\noindent Bipartite graphs determine monomial ideals in the polynomial ring in
two sets of variables $R=K[X_1,\ldots,X_n;Y_1,\ldots,Y_m]$, where
$n$ is the number of the vertices $x_{1},\ldots,x_{n}$ in $V_1$ and
$m$ is the number of the vertices $y_{1},\ldots,y_{m}$ in $V_2$.\\
The edge ideal $I_G$ associated to a bipartite graph $G$ is the
ideal of $R$ that is generated by the monomials of degree two,
$X_iY_j$, on two disjoint sets of variables $X_1, \ldots,
X_n;Y_1,\ldots,Y_m$, such that $\{x_i,y_j\} \in E(G)$, for $1\leq
i\leq n$ and $1\leq j\leq m$.\\
In the following we denote $\mathcal{I} = I_{G}$ when $G=B_{2t}$.

\section{Graded Betti numbers associated to $B_{2t}$}

We are interested in finding bounds for the graded Betti numbers that
appear in the minimal graded resolution of the edge ideal of
$B_{2t}$, in particular  we give upper bounds for them in terms of
the number of the regions of  $B_{2t}$.
\begin{definition}
Let $G$ be a graph on vertex set $V(G)$. We call \emph{induced
subgraph} of $G$ the graph $H \subseteq G$ which has an edge between
any two vertices of it if and only if there is an edge between them
in $G$.
\end{definition}

\begin{proposition} [\cite{J}, 4.1.1 Proposition] \label{3.2}
Let $G$ be a graph. If $H$ is an induced subgraph of $G$ on a subset of
the vertices of $G$, 
then
$$ b_{i_j}(H) \leq  b_{i_j}(G), \quad \forall i,j\,,$$
where $b_{i_j}(H)$ (resp. $ b_{i_j}(G)$) are the graded Betti numbers
associated to $H$ (resp. $G$).
\end{proposition}
\begin{proposition}\label{2.2}
Let  $B_{2t}$ be the bipartite planar graph,  $r=2t$ be the number
of its regions and $\mathcal{I}$ be its edge ideal. Let
$b_{i_j}(B_{2t})$ be the graded Betti numbers in  the minimal
graded resolution of $R/\mathcal{I}$. Then:\\
$1)$ $b_{i_j}(B_{2t}) \leq \sum_{\substack{k+l=i+1\\ k,l\geq 1}}{\frac{3r+12}{4}
\choose k}{\frac{3r+4}{4} \choose l}$,   if $t$ is even;\\
$2)$ $b_{i_j}(B_{2t}) \leq \sum_{\substack{k+l=i+1\\ k,l\geq 1}}{\frac{3r+10}{4}
\choose k}{\frac{3r+6}{4} \choose l}$, if $t$ is odd.
\end{proposition}
\begin{proof} $B_{2t}$ is a   bipartite planar graph on two
disjoint vertex sets $V_1  =\{x_1, \ldots, x_n \}$ and $V_2 =\{y_1,
\ldots, y_m \}$, but it is not complete. Moreover it is an induced
subgraph of the complete bipartite graph on vertex sets
$\overline{V_1} =\{x_1, \ldots, x_{n+1} \}$ and $V_2=\{y_1,
\ldots, y_m \}$.\\
$1)$ If $t$ is even  we have  $|V_1| = n = \frac{3t+4}{2}$ and
$|V_2| = m =\frac{3t+2}{2}$. $B_{2t}$ is an induced subgraph of the
complete bipartite graph $K_{n+1,m}$, where $n+1 =\frac{3t+6}{2}$
and $m=\frac{3t+2}{2}$, that is $V(B_{2t}) \subset V(K_{n+1,m})$ and
$|E(B_{2t})|< |E(K_{n+1,m})|$. Then by Proposition \ref{3.2} we
have: $ b_{i_j}(B_{2t}) \leq b_{i_j}(K_{n+1,m})$, where
$b_{i_j}(K_{n+1,m})$ are the graded Betti numbers of
$R/I(K_{n+1,m})$. By \cite{J}, 5.2.4 Theorem, we have:
$b_{i_j}(K_{n+1,m})= \sum_{\substack{k+l=i+1\\ k,l\geq 1}}{n+1 \choose k}{m
\choose l}.$\\ It follows:
$$b_{i_j}(B_{2t}) \leq \sum_{\substack{k+l=i+1\\ k,l\geq 1}}{\frac{3t+6}{2} \choose k}{\frac{3t+2}{2}
\choose l}, \ \  t=\frac{r}{2}\;.$$
$2)$ If $t$ is odd
 we have  $|V_1| = |V_2| = \frac{3t+3}{2}$. $B_{2t}$ is an induced subgraph of the complete bipartite graph
$K_{n+1,m}$, where $n+1 =\frac{3t+5}{2}$ and    $m=\frac{3t+3}{2}$,
that is $V(B_{2t}) \subset V(K_{n+1,m})$ and $|E(B_{2t})|<
|E(K_{n+1,m})|$.
As before we obtain
$$b_{i_j}(B_{2t}) \leq \sum_{\substack{k+l=i+1\\ k,l\geq 1}}{\frac{3t+5}{2} \choose k}{\frac{3t+3}{2} \choose l}, \ \  t=\frac{r}{2}\;.$$\\[-1.2cm]
\end{proof}
\vspace{3mm}\noindent
It is possible to express the second graded Betti number in degree $3$
of $R/I_G$ in terms of graph theoretical properties for any graph $G$.\\
For a simple graph $G$ there exists the so-called \textit{edge
graph} $L(G)$ of $G$ (\cite{Vill:Real}). It has vertex set equal to
the edge set of $G$ and two vertices of $L(G)$ are adjacent whenever
the corresponding edges of $G$ have one common vertex:
\[
V(L(G))=E(G)=\{f_1,\ldots, f_q \}\]
\[
E(L(G))= \{ (f_i,f_j) \ | \ f_i=\{v_i,v_j\}, \ \  f_j=\{v_j,v_k\}, \ \
i\neq j,\ \  j \neq k \}.
\]
If $G$ is a simple graph on vertices $v_1,\ldots, v_n$ , then the
number of edges of $L(G)$ is given by
$$|E(L(G))|= - |E(G)| + \sum_{i=1}^n\frac{deg^2(v_i)}{2},$$
where $deg(v_i)$ is the number of edges incident with $v_i$.
\begin{theorem} [\cite{El:Vill}] \label{t}
Let $G$  be  a graph and $I_G$ be its edge ideal. If
$$
\ldots \rightarrow R^{c}(-4)\oplus R^{b}(-3)\rightarrow R^{q}(-2)
 \rightarrow R \rightarrow  R/I_G \rightarrow 0
$$
is the minimal graded resolution of $R/I_G$ and $L(G)$ is the edge
graph of  $G$, then
$$b=|E(L(G))|- N_3,$$
where $N_3$ is the number of the triangles of $G$.
\end{theorem}
\noindent In particular, for  $G = B_{2t}$ we can establish the following
\begin{theorem} \label{T2}
Let $B_{2t}$ be the bipartite planar graph, $r=2t$ be the number of
its regions and $\mathcal{I}$ be its edge ideal. If
$$
\ldots \rightarrow R^{c}(-4)\oplus R^{b}(-3)\rightarrow R^{q}(-2)
\rightarrow R \rightarrow R/\mathcal{I} \rightarrow 0
$$
is the minimal graded resolution of $R/\mathcal{I}$, then:
\item{1)} $q= \frac{5}{2}r + 2$;
\item{2)} $b= 6r - 2$.
\end{theorem}
\begin{proof} $1)$ $q=|E(B_{2t})|= |\{ \{v_i,v_{i+1} \} : \ \ 1
\leq i \leq 3t+2, \ \ i\neq t+1, t+2 \}| + |\{ \{v_i,v_{i+t+1}\} : \
\  1 \leq i \leq 2t+2 \}| = (3t+2-2)+ (2t+2)= 5t+2= \frac{5}{2}r +
2$.\\
$2)$ By Theorem \ref{t}, $b=|E(L(B_{2t}))|- N_3$, where $N_3=0$
because the graph is bipartite.
One has:\\
$|E(L(B_{2t}))|= - |E(B_{2t})| + \sum_{i=1}^N\frac{deg^2(v_i)}{2}$,
where $N=3t+3$.\\ We observe that $B_{2t}$ has $N=3(t+1)$ vertices
representable in the plane on three horizontal lines   and on each
line there are $t+1$ vertices. In fact  the representation  in the
plane of $B_{2t}$ is a sequence of squares without chords disposed
in $2$ rows and $t$ columns. It follows that
$$\sum_{i=1}^{3t+3}\frac{deg^2(v_i)}{2}=  4 \left(\frac{2^2}{2}\right) +
2t\left(\frac{3^2}{2}\right) + (t-1)\left(\frac{4^2}{2}\right) =17 t= \frac{17}{2}r\,,$$
where $deg (v_1)= deg(v_{t+1})= deg(v_{2t+3}) = deg(v_{3t+3})=2$,\\
$deg(v_i)=3$ for $2\leq i \leq t$,  $i=t+2,2t+2$   and $2t+4 \leq i\leq 3t+2$,\\
$deg(v_i)=4$ for $t+3 \leq i\leq 2t+1$.\\
Then $\displaystyle{b= |E(L(B_{2t}))|= -\left(\frac{5}{2}r + 2\right) + \frac{17}{2}r=6r - 2\,.}$
\end{proof}

\section{Algebraic topics on minimal vertex covers of $B_{2t}$}

\begin{definition}\rm{
Let $G$ be a graph with vertex set $V(G)$. A subset  $\mathcal{A}
\subset V(G)$ is
called  a  \textit{minimal vertex cover}  for $G$ if:\\
(1) each edge of $G$ is incident with one vertex in $\mathcal{A}$; \\
(2) there is no proper subset  of $\mathcal{A}$ with this property.}
\end{definition}
\noindent If $\mathcal{A}$ satisfies condition (1) only, then
$\mathcal{A}$ is called a \textit{vertex cover} of $G$ and
$\mathcal{A}$ is said to cover all the edges of $G$.\\
The smallest number of vertices in any  minimal vertex cover of $G$
is said \textit{vertex covering number}. We denote it by $\alpha_0(G)$.
\begin{proposition}\label{P}
Let $B_{2t}$ be the bipartite planar graph with $r=2t$ regions and
$t\geq 1$. Then:
$$ {\alpha_0(B_{2t})} = \left\{ \begin{array}{ll}
\frac{3}{4}r + \frac{3}{2}  & \textrm{if $t$ odd}\,\,,\\[2mm]
\frac{3}{4}r + 1  & \textrm{if $t$ even}\,\,.\\
\end{array} \right.
$$
\end{proposition}
\begin{proof} Let $V(B_{2t})=\{v_1, \ldots, v_{3t+3} \}$ and
$E(B_{2t})= \{ \{v_i,v_{i+1} \} \ | \ 1 \leq i \leq 3t+2, \ i\neq t+1,
2t+2, 3t+3 \}  \cup  \{ \{v_i,v_{i+t+1}\} | 1 \leq i \leq 2t+2 \}$.
Hence the representation of $B_{2t}$ in the plane is a sequence of
squares without chords disposed in $2$ rows and $t$ columns.\\[-4mm]

\begin{figure}[htbp]
\begin{center}
   \includegraphics[scale=.65]{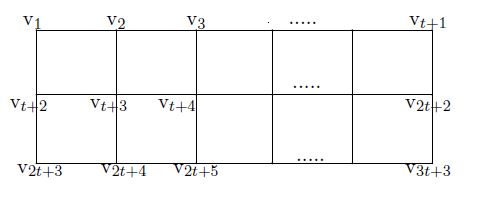}
\end{center}
\end{figure}

\vspace{-6mm}\noindent For $t=1$ and $\alpha_0(B_2)= 3$, it is

\begin{figure}[htbp]
\begin{center}
   \includegraphics[scale=.70]{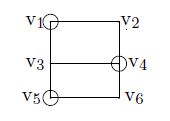}
\end{center}
\end{figure}

\vspace{-6mm}\noindent and $\mathcal{A}(B_2)= \{v_1, v_4, v_5 \} $, $\mathcal{A'}(B_2)=
\{v_2, v_3, v_6 \} $  are minimal vertex covers of $B_{2}$.\\
For $t>1$ a minimal vertex cover of $\alpha_0(B_{2t})$ is given by
adjoining to the minimal vertex cover of $B_2$ one vertex for each
even column and two vertices for each odd column.
Hence\\
$1)$ If $t$ is odd
$$\alpha_0(B_{2t}) =\alpha_0(B_2) + 1\left(\frac{t-1}{2}\right) +
2\left(\frac{t-1}{2}\right)= \frac{3}{2}t + \frac{3}{2}= \frac{3}{4}r +
\frac{3}{2},$$ where $\frac{t-1}{2}$ is the number of the
even (odd) columns in the graph for $t>1$.\\
$2)$ If $t$ is even
$$\alpha_0(B_{2t}) =\alpha_0(B_2) + 1\left(\frac{t}{2}\right) +
2\left(\frac{t}{2}-1\right)= \frac{3}{2}t + 1 =\frac{3}{4}r + 1,$$ where
$\frac{t}{2}$ is the number of the even
columns and $\frac{t}{2} - 1$ is the number
of the odd columns of the graph for $t>1$.
\end{proof}
\begin{proposition}
Let $B_{2t}$ be the bipartite planar graph with $r=2t$ regions and
$t\geq 1$. Then the minimal vertex covers with cardinality
$\alpha_0(B_{2t})$ are:  $$ {\mathcal{A}(B_{2t})} = \left\{
\begin{array}{ll}
V_1 , V_2  & \textrm{if $t$ odd}\,\,,\\
V_2  & \textrm{if $t$ even}\,\,.\\
\end{array} \right.
$$
\end{proposition}
\begin{proof} If $t$ is odd, $\alpha_0(B_{2t})= \frac{3t+3}{2}$,
that is the cardinality of the vertex sets $V_1$ and $V_2$, where
$V_{1}=\{v_{1},v_{3}, \dots, v_{t} \} \cup
\{v_{2+(t+1)},v_{4+(t+1)}, \dots, v_{t+1+(t+1)} \} \cup $ $
\{v_{1+(2t+2)},v_{3+(2t+2)}, \dots, v_{t +(2t+2)} \}$ and
$V_{2}=\{v_{2},v_{4}, \dots, v_{t+1} \} \cup \{v_{1+(t+1)},$ $
v_{3+(t+1)}, \dots, v_{t+(t+1)} \} \cup \{v_{2+(2t+2)},v_{4+(2t+2)},
\dots, v_{t+1 +(2t+2)} \}$ are two
minimal sets of vertices that cover all the edges of $B_{2t}$.\\
If $t$ is even, $\alpha_0(B_{2t})= \frac{3t+2}{2}$, that
is the cardinality of the vertex set  $V_2= \{v_{i} \ | \ i \ {\rm even},\, 1
\leq i \leq 3t+3\}$.  $V_2$ is the only subset of vertices with
cardinality $\alpha_0(B_{2t})$ that covers all the edges of
$B_{2t}$.
\end{proof}
\noindent Now we consider some algebraic aspects linked to the minimal vertex
covers.\\
An immediate consequence of Proposition \ref{P} is the following
\begin{corollary}
Let $\mathcal{I}$ be the edge ideal of $B_{2t}$ with $r=2t$ regions. Then:\\
{$$ {ht(\mathcal{I})} = \left\{ \begin{array}{ll}
\frac{3}{4}r + \frac{3}{2}  & \textrm{if $t$ odd}\,\,,\\[1.5mm]
\frac{3}{4}r + 1  & \textrm{if $t$ even}\,\,.\\
\end{array} \right.
$$}
\end{corollary}
\begin{proof}  It is known that the vertex covering number
$\alpha_0(G)$ is equal to the height of the edge ideal $ht(I_G)$
(\cite{Vill:Real}, 6.1.18). Hence the assertion follows by
Proposition \ref{P}.
\end{proof}
\begin{proposition}
Let $B_{2t}$ be the bipartite planar graph with $r=2t$ regions, $t\geq 1$,
and $\mathcal{I}$ be its edge ideal. Then:$$
{\textrm{dim}(R/\mathcal{I})} = \left\{ \begin{array}{ll}
\frac{3}{4}r + \frac{3}{2}  & \textrm{if $t$ odd}\,\,,\\[1.5mm]
\frac{3}{4}r + 2  & \textrm{if $t$ even}\,\,.\\
\end{array} \right.
$$
\end{proposition}
\begin{proof} Let $R=K[X_1,\ldots,X_n;Y_1,\ldots,Y_m]$ and
$\mathcal{I} \subset R$ be the edge ideal of $B_{2t}$ with
$|V(B_{2t})|=n+m=3t+3$. By \cite{Vill:Real}, (2.1.7), we have $dim(R/
\mathcal{I})=dim(R) -ht(\mathcal{I})$ and
 by \cite{Vill:Real}, (6.1.18), $ht(\mathcal{I})=\alpha_0(B_{2t})$.
Hence $dim(R/
\mathcal{I})=(n+m)-\alpha_0(B_{2t})=3t+3-\alpha_0(B_{2t})$. Then, by
Proposition \ref{P}, it follows:\\
$1)$ $dim(R/ \mathcal{I})= \frac{3}{2}r +3 - (\frac{3}{4}r +
\frac{3}{2})=
\frac{3}{4}r + \frac{3}{2}$, if $t$ is odd,\\
$2)$ $dim(R/ \mathcal{I})= \frac{3}{2}r +3 - (\frac{3}{4}r +
1)=\frac{3}{4}r + 2$, if $t$ is even.
\end{proof}
\begin{proposition}
Let $B_{2t}$ be the bipartite planar graph with $r=2t$ regions and
$\mathcal{I}$ be its edge ideal. Then, for the projective dimension of $R/\mathcal{I}$, we have:\\
$1)$ $\frac{3}{4}r + \frac{3}{2} < pd_R(R/\mathcal{I}) \leq \frac{3}{2}r +3 $, if $t$ is odd;\\
$2)$ $\frac{3}{4}r+1 < pd_R(R/\mathcal{I})\leq \frac{3}{2}r +3$, if
$t$ is even.
\end{proposition}
\begin{proof} For the lower bounds, by \cite{Vill:Real}, Theorem 2.5.14,
one has  $pd_R(R/\mathcal{I}) > ht(\mathcal{I})$. Hence, by
\cite{Vill:Real}, Proposition 6.1.18, it follows $pd_R(\mathcal{I})
> \alpha_0(B_{2t})$. Then the thesis follows by Proposition \ref{P}.\\
For the upper bound, we observe that $B_{2t}$ is an induced
subgraph of the  bipartite complete graph $K_{n+1,m}$, where $n=\frac{3t+4}{2}$  and $m= \frac{3t+2}{2}$ as in Proposition \ref{2.2}.\\
The projective dimension of a graph is affected by some simple
transformations  such as deleting some edges. So, as a consequence
of \cite{J}, Proposition 4.1.3, we have  $pd_R(R/\mathcal{I})  \leq
pd_R(R/I(K_{n+1,m})) =n+1+m-1 = n+m$ (\cite{J}, Proposition 4.2.9),
with $n+m=3t+3= \frac{3}{2}r+3$. Hence: $pd_R(R/\mathcal{I})\leq
\frac{3}{2}r+3$.
\end{proof}
\noindent Finally we recall the one to one correspondence among the minimal
vertex covers of $G$ and minimal primes of $I_G$. In fact, $\wp$ is a minimal prime ideal of
$I_G$ if and only if $\wp=(\mathcal{A})$ for some minimal vertex cover $\mathcal{A}$ of $G$ (\cite{Vill:Real}, 6.1.16). Thus the primary decomposition of the edge ideal of $G$ is given by
$I_G=(\mathcal{A}_1) \cap \cdots \cap (\mathcal{A}_p)$, where
$\mathcal{A}_1, \ldots, \mathcal{A}_p$ are the minimal vertex covers
of $G$. Hence $G$ is an \textit{unmixed} graph if and only if all
the  minimal vertex covers of $G$ have the same cardinality. In
order to study the unmixedness of $B_{2t}$ we recall the following
result:
\begin{proposition} [\cite{Vil}] \label{T}
Let $G$ be a bipartite graph without isolated vertices. Then $G$ is
unmixed if and only if there is a bipartition $V_1=\{x_1, \ldots,
x_m\}$, $V_2=\{y_1, \ldots, y_m\}$ of $G$ such that: \\
1) $\{x_i,y_i \} \in E(G)$ for all $i$;\\
2) if $\{x_i,y_j \}$ and $\{x_j,y_k \}$ are in $E(B_{2t})$, then $\{x_i,y_k \} \in E(B_{2t})$, \,\mbox{$i,j,k$ distinct.}
\end{proposition}
\begin{theorem}
$B_{2t}$ is not unmixed for all $t>0$.
\end{theorem}
\begin{proof} If $t$ is odd, using the characterization of
unmixed bipartite graphs as in the previous proposition, it is
enough to verify that, if $\{x_i,y_j \}, \{x_j,y_k \} \in
E(B_{2t})$, then $\{x_i,y_k \} \notin E(B_{2t})$.\\
Let $V_{1}=\{v_{1},v_{3}, \dots, v_{t} \} \cup
\{v_{2+(t+1)},v_{4+(t+1)}, \dots, v_{t+1+(t+1)} \} \cup
\{v_{1+(2t+2)},$ $v_{3+(2t+2)}, \dots, v_{t +(2t+2)} \}$ and
$V_{2}=\{v_{2},v_{4}, \dots, v_{t+1} \} \cup
\{v_{1+(t+1)},v_{3+(t+1)},  \dots,$ $ v_{t+(t+1)} \} \cup
\{v_{2+(2t+2)},v_{4+(2t+2)}, \dots, v_{t+1 +(2t+2)} \}$ the two
disjoint vertex sets of $B_{2t}$. Replacing with $\{x_1,\ldots,
x_{\frac{3t+3}{2}}\}$  the vertices of $V_1$ and with $\{y_1,\ldots,
y_{\frac{3t+3}{2}}\}$ the vertices of $V_2$, we have $v_1=x_1$,
$v_{1+(t+1)}= y_{\frac{t+1}{2}+1}$, $v_{2+(t+1)}=
x_{\frac{t+1}{2}+1}$, $v_{3+(t+1)}= y_{\frac{t+1}{2}+2}$.\\ Then
$\{x_1, y_{\frac{t+1}{2}+1} \}, \{x_{\frac{t+1}{2}+1},
y_{\frac{t+1}{2}+2}  \} \in E(B_{2t})$, but
$\{x_1,y_{\frac{t+1}{2}+2}  \} \notin E(B_{2t})$.\\
\noindent 2) If $t$ is even, it is sufficient to observe that $V_1$
and $V_2$ are two minimal vertex covers with $|V_1| > |V_2|$.
Hence $B_{2t}$ is not unmixed.
\end{proof}

\section{Perfect matchings of $B_{2t}$}

Let $G$ be a graph. A minimal vertex cover $\mathcal{A}$ of $G$ is
linked to the \textit{set of the independent edges}.
The edges of $G$ that have no common vertices are called \textit{independent edges}.
The \textit{independence number} of a graph $G$, denoted by
$\beta_1(G)$, is the maximum number of its independent edges.
\begin{definition}\rm{
A \textit{matching} of $G$ is a set $\mathcal{M}$ of independent edges.}
\end{definition}
\begin{definition}\rm{
$G$ has a \textit{perfect matching} if it has an even number of
vertices and there is a set of independent edges covering all the
vertices.}
\end{definition}
\noindent This means that there is a pairing off of all the vertices of $G$.
\begin{definition}
A \textit{maximum matching} is a matching $\mathcal{M}$
such that every other matching $\mathcal{M}'$ satisfies
$|\mathcal{M}'| < |\mathcal{M}|$.  In this case $|\mathcal{M}| =
\beta_1(G)$.
\end{definition}
\begin{remark}
Let $\mathcal{M}$ be a maximum matching  and $\mathcal{A}$ a
minimal cover of a graph $G$. Note that each edge of $\mathcal{M}$
must be covered by at least one vertex of $\mathcal{A}$ and each
vertex of $\mathcal{A}$ can cover at most one edge of $\mathcal{M}$.
It follows: $\beta_1(G) \leq \alpha_0(G)$.
\end{remark}
\begin{definition}
A \textit{ perfect matching of K\"{o}nig  type} of a graph $G$ is a
collection $e_1, \ldots, e_g$ of pairwise disjoint edges such that
the union of the vertices in which $e_1, \ldots, e_g$ are incident
is the vertex set of $B_{2t}$ and $g$ is equal to the height of $I_G$.
\end{definition}
\begin{remark}
A graph $G$ satisfies the K\"{o}nig property if the maximum number
of independent edges of $G$ equals the height of $I_G$. Hence a
graph with a perfect matching of K\"{o}nig  type has the K\"{o}nig
property. In \cite{GRV} it is proved that the converse is true for
unmixed graphs.
\end{remark}
\noindent We are interested in analyzing bipartite matching problem, namely in finding a matching with the maximum number of edges.
Clearly, the size of any matching is at most the size of any vertex
cover. This follows from the fact that, given any matching
$\mathcal{M}$, a vertex cover $\mathcal{A}$ must contain at least
one of the vertices of each edge in $\mathcal{M}$. The maximum size
of a matching is at most the minimal cardinality of a vertex cover.

\begin{proposition} [\cite{Vill:Real}, Proposition 6.1.7] \label{1}
For any  bipartite graph $G$,  the size of a maximum matching is
equal to the size of a minimal vertex cover, that is $\beta_1(G) =
\alpha_0(G)$.
\end{proposition}
\begin{theorem}
Let $B_{2t}$ be the bipartite planar graph with $r=2t$ regions,
$t$ odd. Each maximum matching is a perfect matching of cardinality
$\frac{3}{4}r + \frac{3}{2} $.
\end{theorem}
\begin{proof} $B_{2t}$ is a bipartite graph, then, by Proposition
\ref{1}, $\beta_1(B_{2t}) = \alpha_0(B_{2t})$. Hence any vertex of
the minimal vertex cover is incident upon independent edges. Then
$B_{2t}$ has maximum matching with cardinality $\beta_1(B_{2t})=
|\mathcal{M}(B_{2t})|= |V_1|=|V_2|= \frac{3}{4}r + \frac{3}{2}$,
$r=2t.$ Moreover   $B_{2t}$ has an even number of vertices and
$|V_1|=|V_2|$, this means that there is a pairing off of all the
vertices of $B_{2t}$. It follows that each maximum matching  is a perfect matching.
\end{proof}
\begin{theorem}
Let $B_{2t}$ be the bipartite planar graph with $r=2t$ regions,
$t$ odd. $B_{2t}$ has perfect matching of K\"{o}nig type.
\end{theorem}
\begin{proof} $V_{1}=\{v_{1},v_{3}, \dots, v_{t} \} \cup
\{v_{2+(t+1)},v_{4+(t+1)}, \dots, v_{t+1+(t+1)} \} \cup
\{v_{1+(2t+2)},$ $v_{3+(2t+2)}, \dots, v_{t +(2t+2)} \}$ and
$V_{2}=\{v_{2},v_{4}, \dots, v_{t+1} \} \cup \{v_{1+(t+1)},
v_{3+(t+1)}, \dots,$ $ v_{t+(t+1)} \} \cup
\{v_{2+(2t+2)},v_{4+(2t+2)}, \dots, v_{t+1 +(2t+2)} \}$ are  minimal
vertex covers of $B_{2t}$ with cardinality $\alpha_0(B_{2t})=
\frac{3}{4}r + \frac{3}{2}$. Note that
$\beta_1(B_{2t})=\alpha_0(B_{2t})= \frac{3}{4}r + \frac{3}{2}$ and
any vertex of the minimal vertex cover is incident upon independent
edges. Hence, by the geometry of the planar graph $B_{2t}$, we obtain
the following maximum matchings:\\[3mm]
$\bullet \quad \mathcal{M}=\{\{v_{i-1},v_i\} \ | \ i$ {\rm even,} \
$ 2\leq i \leq 3t+3$\}. \vspace{-3mm}

\begin{figure}[htbp]
\begin{center}
   \includegraphics[scale=.65]{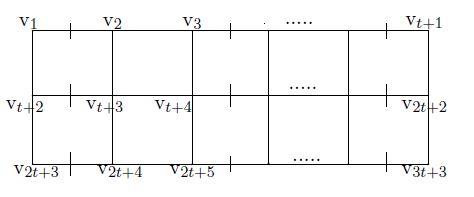}
\end{center}
\end{figure}

\vspace{-2mm}\noindent $\bullet \quad \mathcal{M}= \{\{v_{i-1},v_i\} \ | \ i$ {\rm even,} \ $ 2\leq i \leq t+1 \} \cup
 \{\{v_{i+t},v_{i+2t+1} \} \ | \ 2\leq i \leq t+2\}$. \vspace{-5mm}

\begin{figure}[htbp]
\begin{center}
   \includegraphics[scale=.65]{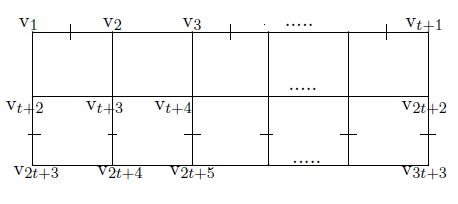}
\end{center}
\end{figure}

\vspace{-2mm}\noindent $\bullet \quad \mathcal{M}= \{\{v_{i},v_{i+t+1}\} \ | \ 1\leq i
\leq t+1 \} \cup \{\{v_{i-1},v_{i} \} \ | \ i$ {\rm even,} \ \mbox{$2t+4 \leq i
\leq 3t+3$\}.} \vspace{-4mm}

\begin{figure}[htbp]
\begin{center}
   \includegraphics[scale=.65]{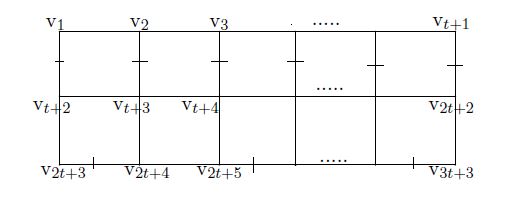}
\end{center}
\end{figure}

\vspace{-6mm}\noindent The other perfect matchings of K\"{o}nig type
are obtained by the previous schemes through different combinations
of the columns in the representation of the graph. In all the cases
$\mathcal{M}$ is a matching such that
$|\mathcal{M}|=\alpha_0(B_{2t})=ht(I(B_{2t}))$ and the union of the
vertices in which the edges of $\mathcal{M}$ are incident coincides
with the vertex set of $B_{2t}$. Hence $B_{2t}$ has perfect
matchings of K\"{o}nig type.
\end{proof}
\noindent Observe that each maximum matching  $\mathcal{M}(B_{2t})$ is a complete
matching from $V_2$ to $V_1$ (being $|V_2|< |V_1|$). This means that
$\mathcal{M}(B_{2t})$ covers each vertex of $V_2$, but not all the
vertices of $V_1$; in fact, $|\mathcal{M}(B_{2t})|= \beta_1(B_{2t})=
|V_2|= \frac{3}{4}r + 1$.

\vspace{.5cm}\par\noindent
{\Large \textbf{Acknowledgements}}\\[2.5mm]
The research that led to the present paper was partially supported by a grant of the group GNSAGA of INdAM, Italy.

\end{document}